\documentclass[12pt]{amsart}
\usepackage[margin=0.9in]{geometry}
\newcommand{\R}{\mathbb{R}}

\newcommand{\N}{\mathbb{N}}

\usepackage{calrsfs}

\renewcommand{\epsilon}{\varepsilon}

\newcommand{\e}{\varepsilon}

\usepackage[english]{babel}
\usepackage[alphabetic]{amsrefs}
\usepackage{amsmath,amssymb,amsfonts,amsthm,enumerate}
\usepackage{url}
\usepackage{graphicx,epstopdf,color}

\usepackage{wrapfig}

\numberwithin{equation}{section}

\newtheorem{theorem}{Theorem}[section]

\newtheorem{corollary}[theorem]{Corollary}

\newtheorem{remark}[theorem]{Remark}

\renewcommand{\leq}{\leqslant}
\renewcommand{\le}{\leqslant}

\renewcommand{\ge}{\geqslant}

\begin{document}

\title[Degenerate double-well potentials]{Density estimates for degenerate double-well potentials}

\author[Serena Dipierro, Alberto Farina and Enrico Valdinoci]{
Serena Dipierro${}^{(1)}$
\and
Alberto Farina${}^{(2)}$
\and
Enrico Valdinoci${}^{(1,3, 4)}$
}

\subjclass[2010]{35J61, 35J62, 35J70, 35J75
35J15, 35J20, 82B26.}
\keywords{Allen-Cahn equation, Cahn-Hilliard equation, phase coexistence models,
non-degeneracy assumptions.}

\maketitle

{\scriptsize \begin{center} (1) -- Dipartimento di Matematica ``Federigo Enriques''\\
Universit\`a degli studi di Milano\\
Via Saldini 50, I-20133 Milano (Italy)\\
\end{center} 
\scriptsize \begin{center} (2) --
LAMFA,
CNRS UMR 7352\\
Facult\'e des Sciences\\
Universit\'e de Picardie Jules Verne\\
33 rue Saint Leu,
80039 Amiens CEDEX 1 (France)\end{center}
\scriptsize \begin{center} (3) -- School of Mathematics and Statistics  \\
University of Melbourne\\ Grattan Street, 
Parkville, VIC-3010 Melbourne (Australia)\\
\end{center}
\scriptsize \begin{center} (4) --
Istituto di Matematica Applicata e Tecnologie Informatiche\\
Via Ferrata 1, 27100 Pavia (Italy)\\ \end{center}
\bigskip

\begin{center}
E-mail addresses:
{\tt serena.dipierro@unimi.it},
{\tt alberto.farina@u-picardie.fr},
{\tt enrico@mat.uniroma3.it}
\end{center}
}

\begin{abstract}
We consider a general energy functional for phase coexistence models,
which comprises the case of Banach norms in the gradient term
plus a double-well potential.

We establish density estimates for $Q$-minima. Namely,
the state parameters close to both phases are proved to occupy
a considerable portion of the ambient space.

{F}rom this, we obtain the uniform convergence of the level sets to the limit
interface in the sense of Hausdorff distance.

The main novelty of these results lies in the fact that we do not assume
the double-well potential to be non-degenerate in the vicinity of the minima.

As far as we know, these types of density results for degenerate potentials
are new even for minimizers and even in the case of semilinear equations, but
our approach can comprise at the same time
quasilinear equations, $Q$-minima and 
general energy functionals.
\end{abstract}

\section{Introduction}

Phase coexistence models study
the separation interfaces between regions corresponding to different values
of a suitable state parameter. A typical model, introduced by
J. D. van der Waals~\cite{MR523642} and developed, under various perspectives,
by S. M. Allen, 
J. W. Cahn and
J. E. Hilliard~\cites{CH, AC}, and (in the vectorial setting) by
V. L. Ginzburg, L. D. Landau and L. P. Pitaevski\u\i~\cites{MR0105929, MR0237287},
considers a domain~$\Omega\subseteq\R^n$ and a state parameter~$u:\Omega\to[-1,1]$.
In such model, the ``pure phases'' correspond to the values of the state parameters~$1$
and~$-1$ and the phase separation is driven by the minimization of
an energy functional.

The typical energy functional taken into account is the superposition of
a potential energy induced by a ``double-well'' function~$W$,
which tries to force the system into the pure phases, and an interaction energy
(e.g. a gradient penalization)
which avoids the production of unnecessary interfaces. More precisely,
the potential energy is often taken of the form
\begin{equation*} {\mathcal{P}}_\Omega(u):=\int_\Omega W(u(x))\,dx\end{equation*}
with~$W\ge0$ and~$W(\tau)>0=W(-1)=W(1)$ for any~$\tau\in(-1,1)$, 
and a natural candidate for the interaction energy is the Dirichlet form
$$ {\mathcal{I}}_\Omega(u):=\frac12\int_\Omega |\nabla u(x)|^2\,dx.$$
In this setting, the total energy becomes
\begin{equation}\label{IP} {\mathcal{I}}_\Omega(u)+{\mathcal{P}}_\Omega(u)
=\frac12\int_\Omega |\nabla u(x)|^2\,dx
+\int_\Omega W(u(x))\,dx\end{equation}
and the critical points are solutions of
\begin{equation} \label{EQG}\Delta u(x)=W'(u(x)) \;{\mbox{ for any }}x\in\Omega.\end{equation}
Probably, the most commonly studied case is that in which~$W(\tau)=\frac{|1-\tau^2|^2}4$:
in this case, formula~\eqref{EQG} reduces to the so-called
Allen-Cahn equation
\begin{equation*} \Delta u(x)+u(x)-u^3(x)=0 \;{\mbox{ for any }}x\in\Omega.\end{equation*}
A natural problem in phase separation models is then to ``describe
the picture seen from afar'': namely, one can expect that, at a large scale,
the two phases tend to ``separate'' one from the other, with ``the least possible interface''.
\medskip

Till now, two main approaches have been adopted to rigorously describe
this phase separation. The first method relies on the theory of $\Gamma$-convergence,
and aims at identifying a limit functional and preserving the notion of minimizers.
The second method is based on density estimates, that is on measuring
which portion of the domain the two phases occupy, and gives as a byproduct
the convergence of the level sets in the Hausdorff distance.

Very roughly speaking, both the methods of~$\Gamma$-convergence
and density estimates consider the scaled problem
obtained by a spatial dilation of a minimizer. Namely, if~$u$ is a minimizer
of the energy functional and~$\e\in(0,1)$ is a small parameter, one
considers
$$ u_\e(x):=u\left(\frac{x}{\e}\right).$$
One of the results of the $\Gamma$-convergence theory,
as established in~\cite{MR0445362}
and greatly extended in~\cites{MR1036589, MR1097327},
is that, as~$\e\searrow0$, the function~$u_\e$ converges (up to a subsequence)
in~$L^1_{\rm loc}(\R^n)$
to a function~$u_0$ which takes values only in~$1$ and~$-1$ (i.e.
the state parameters of~$u_0$
are only pure phases); remarkably, if one defines~$E:=\{ u_0=1\}$,
then~$\partial E$ is a minimal surface, namely~$E$ is a local minimizer
of the perimeter functional.
In addition, a suitable rescaling of the original energy functional
possesses appropriate convergence properties to the perimeter functional,
and this convergence is compatible with energy minimization
(see also~\cite{MR1968440} for a general introduction to this topic).
\medskip

The second approach, based on density estimates,
has been introduced in~\cite{CC},
and several extensions also to inhomogeneous,
singular and degenerate equations were performed in~\cites{MR2099113, MR2126143, MR2139200}.
A general approach to density estimates was also presented in~\cite{FV},
and the vector-valued case
has been treated in~\cite{MR3393249}.
See also the recent monograph~\cite{ALI} for more exhaustive discussions
on density estimates for phase transitions.\medskip

The basic idea of the density estimates
is to consider a point, say the origin, which belongs to the interface, say
for concreteness~$u(0)=0$, and measure the proportion of the sets~$\{u>1/2\}$
and~$\{u<-1/2\}$ in a large ball.
In this configuration, one aims at showing that the Lebesgue measure
of~$B_{r}\cap\{ u\in[-1/2,1/2]\}$ is bounded from above by~$O(r^{n-1})$,
while the Lebesgue measures
of~$B_{r}\cap\{ u>1/2\}$ and~$B_{r}\cap\{ u<-1/2\}$
are bounded from below by~$O(r^n)$ for large~$r$. That is,
state parameters close to
the pure phases occupy a considerable portion of the space on a large scale,
while the interface becomes relatively negligible in measure.
\medskip

A consequence of these density estimates is also that the level sets of~$u_\e$
converge locally in the Hausdorff distance to~$\partial E$: more precisely,
for any~$\delta$, $R>0$, there exists~$\e_0(\delta,R)$ such that
if~$\e\in \big(0,\,\e_0(\delta,R)\big)$ then
$$ \{ u_\e\in(-1/2,1/2)\}\cap B_R\subseteq \bigcup_{x\in \partial E} B_\delta(x).$$

It is clear that the $\Gamma$-convergence and the density estimates
approaches share a common interest
in the large-scale phase separation regimes, but address this problem
with different methods and obtaining different results. In addition,
an important technical difference between these two approaches arises
in the main assumptions on the double-well potential.
Indeed, while the assumptions on the double-well potential
for the $\Gamma$-convergence results are very mild and do not
involve any quantitative hypothesis, the density estimates
usually assume a non-degenerate growth from the minima of the potential.
Roughly speaking, the potential~$W$ is supposed to have
a ``sufficiently strong'' growth from~$\pm1$.
For instance, in the setting of~\eqref{IP}, one can consider the case in which~$W(\tau)\sim |1-\tau^2|^m$,
with~$m\in(0,2]$ (or even~$W(\tau)=\chi_{(-1,1)}(\tau)$), but the case~$m>2$
has never been considered in the literature to the best of our knowledge.
The analytic counterpart of this range of~$m$
is that the known density estimates 
till now rely on the non-degeneracy condition~$W''(\pm1)\ne0$.\medskip

\begin{wrapfigure}{R}{9cm}
    \centering
    \includegraphics[width=9cm]{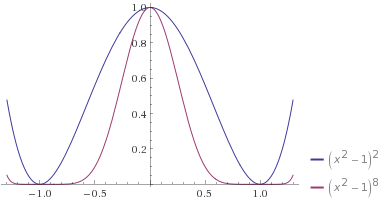}
    \caption{\it {{Non-degenerate versus degenerate potentials.}}}
    \label{WIK1}
\end{wrapfigure}

The case~$W''(\pm1)=0$ is conceptually more complicated since, from
the physical point of view, if~$W$ is ``too flat'', then the system
might develop ``approximative or intermediate phases''. The extreme case of this would
be, for instance, the case in which~$W(\tau)>0=W(\sigma)$ for any~$\tau
\in(-a,a)$ and any~$\sigma\in[-1,-a]\cup[a,1]$,
with~$a\in(0,1)$: in this situation, minimizers would prefer the state parameters~$\pm a$
rather than~$\pm1$ and the state parameters with value in~$(a,1]$ (as well as in~$[-1,-a)$)
become empty, thus violating the density estimates.
\medskip

Interestingly, non-degeneracy assumptions on the potentials also play
an important role in the asymptotic theory of critical points
in terms of the limiting varifold, see in particular Assumption~A on page~51
of~\cite{MR1803974} (moreover, a different
growth from the potential wells can also create different
phenomena in other settings, see e.g.~\cite{MR3385194}
for the case of singular potentials).\medskip

The goal of this paper is to provide density estimates for degenerate potentials,
namely when~$W''(\pm1)=0$, under weaker growth assumptions.
The method of the
proof is very general, it also applies to quasiminima, and we can deal with 
an energy functional, which is possibly non-quadratic
in the gradient term (and so corresponding
also to quasilinear PDEs), that we now introduce
in details.\medskip

Given~$n\in\N$, $n\ge1$, and a bounded, open set~$\Omega\subset\R^n$,
we consider an energy functional of the form
\begin{equation*}
{\mathcal{E}}_\Omega(u):=\int_\Omega 
E\big( x,\,u(x),\,\nabla u(x)\big)\,dx
.\end{equation*}
We suppose that,
for any~$x\in\R^n$, $\tau\in\R$ and~$\xi\in\R^n$, 
it holds that
\begin{eqnarray}
\label{SYT:ST1}
&&E( x,\tau,\xi)\ge \lambda \Big( |\xi|^p + |\tau+1|^{{ m }}
\chi_{[-\infty,\theta]}(\tau)\Big),
\\ \label{SYT:ST2}
{\mbox{and }}&&
E( x,\tau,\xi)\le \Lambda \Big( |\xi|^p + |\tau+1|^{{ m }}
\Big),
\end{eqnarray}
for some~$\theta\in(-1,1)$, 
$\lambda\in(0,1]$,
$\Lambda\in[1,+\infty)$, $p\in (1,+\infty)$ and~${{ m }}\in(p,+\infty)$.
\medskip

Roughly speaking,
conditions~\eqref{SYT:ST1}
and~\eqref{SYT:ST2} mean that the energy density under consideration
is bounded from both sides by a ``gradient term of $L^p$-type''
and, in the vicinity of a fixed phase (say, the phase corresponding to~$-1$),
by ``a well which grows like a power~$m$''
(when~$p=2$, the degeneracy of the potential reflects\footnote{As
a technical observation, we observe that
we will focus our analysis
on each separate well of the potential. In this way,
our result in Theorem~\ref{DENS} also applies to double well-potentials in which
only one of the two wells possesses a growth condition
as in~\eqref{CONDIZIONE}.}
into the fact that~$m>2$,
see Figure~\ref{WIK1}).\medskip

Given~$Q\in[1,+\infty)$, we say that~$u\in W^{1,p}(\Omega)$
is a $Q$-minimum in~$\Omega$ 
if, for any open and bounded set~$\Omega'\Subset\Omega$,
we have that~${\mathcal{E}}_{\Omega'}(u)<+\infty$
and 
\begin{equation}\label{QMI}
{\mathcal{E}}_{\Omega'}(u)\le Q\, {\mathcal{E}}_{\Omega'}(v),
\end{equation}
for any~$v\in W^{1,p}(\Omega')$
with~$u-v\in W^{1,p}_0(\Omega')$.\medskip

Of course, when~$Q=1$, the notion in~\eqref{QMI} boils
down to local minimality.
In the framework of~\eqref{QMI}
(and denoting, as customary,
the $n$-dimensional Lebesgue measure by~${\mathcal{L}}^n$), our main result is the following:

\begin{theorem}\label{DENS}
Let~$u\in W^{1,\infty}(\Omega, \,[-1,1])$ be a $Q$-minimum in~$\Omega$,
with
\begin{equation}\label{INIZIO}
{\mathcal{L}}^n\big(B_{r_o}\cap\{ u>\theta\}\big)\ge c_o
\end{equation}
for some~$r_o$, $c_o>0$, with~$B_{r_o}\subset\Omega$.

Assume that
\begin{equation}\label{CONDIZIONE}
\frac{pm}{m-p} > n.
\end{equation}
Then, there exist~$\hat r$, $\hat c>0$ such that for any~$r\ge\hat r$ for which~$B_r\subset\Omega$
it holds that
\begin{equation}\label{FINE}
\int_{B_{r}\cap\{ u\le\theta\}} |u(x)+1|^m\,dx+
{\mathcal{L}}^n\big(B_{r}\cap\{ u>\theta\}\big)\ge \hat c\,r^n
.\end{equation}
\end{theorem}

\begin{remark} \label{P2345sa}{\rm The non-degenerate case of
Theorem~\ref{DENS}
corresponds to the ranges~$p=2$ and~$m\in(0,2]$
or, more generally,~$p>1$ and~${{ m }}\in(0,p]$,
and it has been dealt with in~\cite{FV}.}
\end{remark}

\begin{remark} {\rm To the best of our knowledge,
the degenerate potential case of 
Theorem~\ref{DENS} is new even when~$p=2$
and for the case of minimizers (instead of quasiminima).
It is also new even for energy functionals of the type
$$ \frac12 \int_\Omega|\nabla u(x)|^2\,dx+\int_\Omega |1-u^2(x)|^m\,dx,$$
with
$$ m\in \left\{\begin{matrix}
(2,+\infty) & {\mbox{ if }}n\in\{1,2\},\\
\, \\
\left(2,\displaystyle\frac{2n}{n-2}\right)&{\mbox{ if }}n\ge3.
\end{matrix}\right.$$}
\end{remark}

\begin{remark}{\rm Condition~\eqref{CONDIZIONE} is always satisfied when~$n=1$.
}\end{remark}

\begin{remark}\label{93243504567}{\rm Condition~\eqref{CONDIZIONE} is always satisfied when~$n\le p$.
}\end{remark}

\begin{remark}{\rm 
It is interesting to observe that
condition~\eqref{CONDIZIONE} is a requirement
relating the critical Sobolev exponent of~$W^{1,p}(\R^m)$ with the dimension
of the ambient space.

Interestingly (assuming~$m$ and~$n$ larger than~$p$
and recalling Remark~\ref{P2345sa} for the case~$m\le p$
and Remark~\ref{93243504567} for the case~$n\le p$), 
we point out that
condition~\eqref{CONDIZIONE} can
also be written as
$$ \frac{np}{n-p}>m,$$
which means that the Sobolev exponent of~$W^{1,p}(\R^n)$
is larger than~$m$ (i.e., the exponent~$m$ in~\eqref{CONDIZIONE}
is assumed to be ``subcritical'').
}\end{remark}

\begin{remark}{\rm 
We think that it is an interesting open problem
to check whether or not
condition~\eqref{CONDIZIONE} is sharp.
}\end{remark}

\begin{remark}{\rm 
We also think that it is a very interesting open problem to establish whether or not
the asymptotic theory of critical points
in terms of the limiting varifold
established in~\cite{MR1803974} remains valid in the case of degenerate potentials.
}\end{remark}

{F}rom Theorem~\ref{DENS}, we obtain the following general density estimate:

\begin{corollary}\label{CORO1}
Assume that for any~$\mu\in(-1,1)$ there exists~$\lambda_\mu>0$
such that, for any~$x\in\R^n$, $\tau\in\R$ and~$\xi\in\R^n$, 
\begin{equation}\label{N}
E(x,\tau,\xi)\ge \lambda_\mu(\tau+1)^m\chi_{(-\infty,\mu]}(\tau).
\end{equation}
Let~$u\in W^{1,\infty}(\Omega, \,[-1,1])$ be a $Q$-minimum in~$\Omega$.
Then, for any~$r>1$ for which~$B_r\subset\Omega$
it holds that
\begin{equation}\label{n-1}
{\mathcal{E}}_{B_r}(u)\le C r^{n-1},
\end{equation}
for some~$C>0$, possibly depending on~$\lambda$, $\Lambda$,
$\theta$, $n$, ${{ m }}$, $p$, $Q$
and~$\| u\|_{W^{1,\infty}(\Omega)}$.

Moreover, if
\begin{equation}\label{0q9fg09weuty3049yt4}
{\mathcal{L}}^n\big(B_{r_o}\cap\{ u>\vartheta\}\big)\ge c_o
\end{equation}
for some~$\vartheta\in(-1,1)$ and~$r_o$, $c_o>0$, with~$B_{r_o}\subset\Omega$, and
\begin{equation*}
\frac{pm}{m-p} > n,
\end{equation*}
then, for any~$\Theta\in(-1,1)$
there exist~$\hat r$, $\hat c>0$,
possibly depending on~$\lambda$, $\Lambda$,
$\theta$, $\vartheta$, $\Theta$, $n$, ${{ m }}$, $p$, $Q$, $r_o$, $c_o$
and~$\| u\|_{W^{1,\infty}(\Omega)}$,
such that for any~$r\ge\hat r$ for which~$B_r\subset\Omega$
it holds that
\begin{equation}\label{LAS}
{\mathcal{L}}^n\big(B_{r}\cap\{ u>\Theta\}\big)\ge \hat c\,r^n
.\end{equation}
\end{corollary}

As a consequence of Corollary~\ref{CORO1} we also have a convergence
result in the Hausdorff distance:

\begin{corollary}\label{CORO2}
Assume that
\begin{eqnarray}
\label{SYT:ST1.1}
&&E( x,\tau,\xi)\ge \lambda_\star\, \big( |\xi|^p + |1-\tau^2|^{{ m }}
\big),
\\ \label{SYT:ST2.1}
{\mbox{and }}&&
E( x,\tau,\xi)\le \Lambda_\star \big( |\xi|^p + |1-\tau^2|^{{ m }}\big),
\end{eqnarray}
for some
$\lambda_\star\in(0,1]$,
$\Lambda_\star\in[1,+\infty)$, $p\in (1,+\infty)$ and~${{ m }}\in(p,+\infty)$. Suppose that
$$ \frac{pm}{m-p}>n.$$

Consider an infinitesimal sequence of~$\e$'s and let~$u_\e\in W^{1,\infty}
(\Omega, [-1,1])$
be a sequence of
$Q$-minima in~$\Omega$ for the functional
$$ {\mathcal{E}}_\Omega^\e (v):=\frac1\e\,\int_\Omega E\big(x, \,v(x),\, \e\nabla v(x)\big)\,dx,$$
with
\begin{equation}
\label{SUPF11} \sup_{\e} \e\|\nabla u_\e\|_{L^\infty(\Omega)}<+\infty\end{equation}
and
\begin{equation}
\label{SUPF} \sup_{\e}  {\mathcal{E}}_\Omega^\e (u_\e)<+\infty.\end{equation}
Then, there exists~$E\subseteq\Omega$ such that, up to a subsequence,
$u_\e$ converges to~$\chi_E-\chi_{\R^n\setminus E}$ in~$L^1_{\rm loc}(\Omega)$.

Furthermore, given any~$\Theta\in(0,1)$, the set~$\{|u_\e|\leq\Theta\}$
converges to~$\partial E$ in the Hausdorff distance locally uniformly in~$\Omega$, namely
for any~$\delta$, $R>0$, for which~$B_R\Subset\Omega$,
there exists~$\e_0(\delta,R,\Theta)$ such that
if~$\e\in \big(0,\,\e_0(\delta,R,\Theta)\big)$ then
\begin{equation} \label{89:098:98rt}
\{ |u_\e|\le\Theta\}\cap B_R\subseteq \bigcup_{x\in \partial E} B_\delta(x).\end{equation}
\end{corollary}

We notice that conditions~\eqref{SYT:ST1.1}
and~\eqref{SYT:ST2.1} are just a rephrasing of~\eqref{SYT:ST1}
and~\eqref{SYT:ST2} that take into account the growth from
both the wells of the potential that induces the two phases~$1$ and~$-1$.\medskip

Next section is devoted to the proof of the main result
of Theorem~\ref{DENS}. Then, the simple proofs
of Corollaries~\ref{CORO1}
and~\ref{CORO2} will be
given in Sections~\ref{SCORO1} and~\ref{SCORO2}, respectively.

\section{Proof of Theorem~\ref{DENS}}

We fix an additional parameter~$T\in(0,+\infty)$, which will be suitably
chosen, possibly also in dependence of~$\lambda$, $\Lambda$,
$\theta$, $n$, ${{ m }}$, $p$, $Q$, $r_o$ and~$c_o$ (recall the
structural constants in~\eqref{SYT:ST1},
\eqref{SYT:ST2} and~\eqref{INIZIO}).

The proof of Theorem~\ref{DENS} is based on an iteration argument.
To this end, for any~$k\in\N$ and~$x\in B_{(k+1)T}$, we set
\begin{equation} \label{vk:owue}
v_k(x):=\frac{2}{\big( 1
+(k+1)T-|x|
\big)^{\frac{p}{m-p}}}-1.\end{equation}
{F}rom now on, $C>0$ will be a constant, which possibly 
depends
on the structural constants~$\lambda$, $\Lambda$,
$\theta$, $n$, ${{ m }}$ and~$p$ in~\eqref{SYT:ST1}
and~\eqref{SYT:ST2},
and also on the fixed~$Q$, $r_o$ and~$c_o$ of
Theorem~\ref{DENS},
but it is independent of~$T$,
and can freely vary from line to line (later on, $T$ will be supposed to be sufficiently
large with respect to~$C$).
In this way, we have that, for any~$x\in B_{(k+1)T}$,
\begin{equation}\label{2.3A}
|\nabla v_k(x)|^p\;=\;C\,\left(
\frac{1}{\big( 1
+(k+1)T-|x|
\big)^{\frac{m}{m-p}} }\right)^p
\;=\; C\,\big( v_k(x)+1\big)^m.
\end{equation}
{F}rom this and~\eqref{SYT:ST2}, 
letting
\begin{equation}\label{2.1BIS}\Omega_k:=B_{(k+1)T}\cap \{u>v_k\},
\end{equation}
we obtain that
\begin{equation}\label{qdwtcg2euf0345}
\begin{split}
{\mathcal{E}}_{\Omega_k}(v_k)\,&=\int_{\Omega_k} 
E\big( x,\,v_k(x),\,\nabla v_k(x)\big)\,dx\\
&\le C \int_{\Omega_k} 
|\nabla v_k(x)|^p + |v_k(x)+1|^m\,dx\\
&\le C \int_{\Omega_k} |v_k(x)+1|^m\,dx
.\end{split}\end{equation}
Notice now that
\begin{equation}\label{89:0234hfhfhfh}
{\mbox{$v_k=1\ge u$ on~$\partial B_{(k+1)T}$,}}
\end{equation}
and therefore, by the definition 
of~$\Omega_k$ in~\eqref{2.1BIS}
and that
of~$v_k$ in~\eqref{vk:owue},
we have that~$v_k=u$ on~$\partial \Omega_k$.
Hence, the $Q$-minimality of~$u$
(recall~\eqref{QMI}) implies that
\begin{equation*}
{\mathcal{E}}_{\Omega_k}(u)\le Q\, {\mathcal{E}}_{\Omega_k}(v_k).
\end{equation*}
Using this, \eqref{SYT:ST1} and~\eqref{qdwtcg2euf0345},
it follows that
\begin{equation}\label{283458757hfhfhfhfhfh8234}
\int_{\Omega_k} 
|\nabla u(x)|^p + |u(x)+1|^m\,\chi_{(-\infty,\theta]}(u(x))\,dx
\le C \int_{\Omega_k} |v_k(x)+1|^m\,dx
.
\end{equation}
On the other hand, by Young's Inequality
(with exponents~$p$ and~$\frac{p}{p-1}$)
and Coarea Formula, for any Lipschitz function~$w$ we have that
\begin{equation}\label{uIAKLAP}
\begin{split}
&C\,\int_{\Omega_k} 
|\nabla w(x)|^p + |w(x)+1|^m\,\chi_{(-\infty,\theta]}(w(x))\,dx\\ \ge\;&
\int_{\Omega_k} 
|\nabla w(x)||w(x)+1|^{\frac{m(p-1)}p}\,\chi_{(-\infty,\theta]}(w(x))\,dx\\
=\;&
\int_\R\left[\int_{\{w=t\}}
\chi_{\Omega_k}(x)\, |w(x)+1|^{\frac{m(p-1)}p}\,\chi_{(-\infty,\theta]}(w(x))\,d{\mathcal{H}}^{n-1}(x)
\right]\,dt\\=\;&
\int_{-\infty}^{\theta}
{\mathcal{H}}^{n-1}\big(\Omega_k\cap\{w=t\}\big)\, |t+1|^{\frac{m(p-1)}p}
\,dt.
\end{split}\end{equation}
Moreover, from~\eqref{vk:owue},
\begin{equation*} 
\sup_{x\in B_{kT}} v_k(x)\le\frac{2}{\big( 1
+T
\big)^{\frac{p}{m-p}}}-1.
\end{equation*}
In particular, choosing~$T$ conveniently large, it holds that
\begin{equation} \label{9er7o2freyurtuwh}
v_k<
\frac{\theta-1}{2} \quad{\mbox{ in }}B_{kT}. \end{equation}
As a consequence,
we have that
\begin{equation*} 
B_{kT}\cap\{u>\theta\}\;\subseteq \;
B_{kT}\cap\{ u>t>v_k\}\end{equation*}
for any~$t\in\left[ \frac{\theta-1}{2},\,\theta\right]$.

In particular, by~\eqref{2.1BIS},
\begin{equation} \label{0987ohg763yedj923484}
B_{kT}\cap\{u>\theta\}\;\subseteq \;
B_{(k+1)T}\cap\{ u>t>v_k\}\;=\;\Omega_k
\cap\{ u>t>v_k\}.
\end{equation}
for any~$t\in\left[ \frac{\theta-1}{2},\,\theta\right]$.

In addition, recalling~\eqref{89:0234hfhfhfh},
we have that~$(\partial B_{(k+1)T})\cap\{ u>t>v_k\}
\subseteq \{v_k=1\ge u\}\cap\{u>v_k\}$
and so
$$ (\partial B_{(k+1)T})\cap\{ u>t>v_k\}=\varnothing.$$
Consequently,
\begin{equation}\label{uyeugdhhdhhfd43hf}\begin{split}
&\partial
\big( \Omega_k\cap\{ u>t>v_k\}\big)
=\partial
\big( B_{(k+1)T}\cap\{ u>t\}\cap\{t>v_k\}\big)\\&\qquad
\subseteq
\big(
B_{(k+1)T}\cap(\partial\{ u>t\})\cap\{t>v_k\}
\big)\cup
\big(
B_{(k+1)T}\cap\{ u>t\}\cap(\partial\{t>v_k\})
\big)\\
&\qquad=
\big(
B_{(k+1)T}\cap\{u=t>v_k\}
\big)\cup
\big(
B_{(k+1)T}\cap\{ u>t=v_k\}
\big)\\
&\qquad=
\big(
\Omega_k\cap\{u=t>v_k\}
\big)\cup
\big(
\Omega_k\cap\{ u>t=v_k\}
\big)\\
&\qquad\subseteq
\big(
\Omega_k\cap\{u=t\}
\big)\cup
\big(
\Omega_k\cap\{ v_k=t\}
\big)
.
\end{split}
\end{equation}
Now we observe that, by Isoperimetric Inequality,
\begin{equation*}
{\mathcal{L}}^n\big( \Omega_k\cap\{ u>t>v_k\}\big)^{\frac{n-1}n}
\le C {\mathcal{H}}^{n-1}\Big(\partial
\big( \Omega_k\cap\{ u>t>v_k\}\big)\Big).
\end{equation*}
This and~\eqref{uyeugdhhdhhfd43hf} imply that
\begin{equation}\label{uIAKLAP2}
{\mathcal{L}}^n\big( \Omega_k\cap\{ u>t>v_k\}\big)^{\frac{n-1}n}
\le C \Big( {\mathcal{H}}^{n-1}
\big(
\Omega_k\cap\{u=t\}
\big)+{\mathcal{H}}^{n-1}
\big(
\Omega_k\cap\{ v_k=t\}
\big)\Big).
\end{equation}
Accordingly, using
first~\eqref{uIAKLAP2},
and then~\eqref{uIAKLAP}, 
with~$w:=u$ and~$w:=v_k$, we conclude that
\begin{equation*}\begin{split}
&
\int_{(\theta-1)/2}^\theta
{\mathcal{L}}^n\big( \Omega_k\cap\{ u>t>v_k\}\big)^{\frac{n-1}n}
\,dt\\
\le\;\,&C\,
\int_{(\theta-1)/2}^\theta
{\mathcal{L}}^n\big( \Omega_k\cap\{ u>t>v_k\}\big)^{\frac{n-1}n}
\, |t+1|^{\frac{m(p-1)}p}
\,dt
\\
\le\;\,&C\,
\int_{-\infty}^\theta
{\mathcal{L}}^n\big( \Omega_k\cap\{ u>t>v_k\}\big)^{\frac{n-1}n}
\, |t+1|^{\frac{m(p-1)}p}
\,dt
\\ 
\le\;\,&C\,\left(
\int_{-\infty}^\theta
{\mathcal{H}}^{n-1}
\big(\Omega_k\cap\{u=t\}\big)\, |t+1|^{\frac{m(p-1)}p}
\,dt+
\int_{-\infty}^\theta
{\mathcal{H}}^{n-1}
\big(\Omega_k\cap\{v_k=t\}\big)\, |t+1|^{\frac{m(p-1)}p}
\,dt\right)
\\ 
\le\;\,
& C\,\left(\int_{\Omega_k} 
|\nabla u(x)|^p + |u(x)+1|^m\,\chi_{(-\infty,\theta]}(u(x))\,dx
+
\int_{\Omega_k} 
|\nabla v_k(x)|^p + |v_k(x)+1|^m\,\chi_{(-\infty,\theta]}(v_k(x))\,dx\right)
.\end{split}\end{equation*}
This and~\eqref{0987ohg763yedj923484} imply that
\begin{eqnarray*}
&&C\,\left(\int_{\Omega_k} 
|\nabla u(x)|^p + |u(x)+1|^m\,\chi_{(-\infty,\theta]}(u(x))\,dx
+
\int_{\Omega_k} 
|\nabla v_k(x)|^p + |v_k(x)+1|^m\,\chi_{(-\infty,\theta]}(v_k(x))\,dx\right)\\
&&\qquad\qquad\qquad
\ge
{\mathcal{L}}^n\big(B_{kT}\cap\{u>\theta\}\big)^{\frac{n-1}n}.
\end{eqnarray*}
Consequently, recalling~\eqref{2.3A}
and~\eqref{283458757hfhfhfhfhfh8234}, we find that
\begin{equation}\label{LAKKA:01284}
{\mathcal{L}}^n\big(B_{kT}\cap\{u>\theta\}\big)^{\frac{n-1}n}\le
C \int_{\Omega_k} |v_k(x)+1|^m\,dx.
\end{equation}
Furthermore, exploiting~\eqref{283458757hfhfhfhfhfh8234},
we see that
\begin{eqnarray*}
&& \int_{\Omega_k\cap\{u\leq\theta\}} |u(x)+1|^m\,dx
=
\int_{\Omega_k} 
|u(x)+1|^m\,\chi_{(-\infty,\theta]}(u(x))\,dx
\\ &&\qquad\le
\int_{\Omega_k} 
|\nabla u(x)|^p + |u(x)+1|^m\,\chi_{(-\infty,\theta]}(u(x))\,dx
\le C \int_{\Omega_k} |v_k(x)+1|^m\,dx
.\end{eqnarray*}
As a consequence of this, and
recalling~\eqref{2.1BIS} and~\eqref{9er7o2freyurtuwh},
we have 
\begin{equation}\label{9qw8e7ryt32938yuwhjsdhfg}
\begin{split}&
\int_{B_{kT}\cap \{u\le\theta\}} |u(x)+1|^m\,dx\\ =\;\,&
\int_{B_{kT}\cap \{\theta\ge u> v_k\}} |u(x)+1|^m\,dx+
\int_{B_{kT}\cap \{u\le\min\{v_k,\theta\}\}} |u(x)+1|^m\,dx\\
\le\;\,&
\int_{B_{kT}\cap \{\theta\ge u>v_k\}} |u(x)+1|^m\,dx+
\int_{B_{kT}\cap \{u\le v_k\}} (u(x)+1)^m\,dx\\
\le\;\,&
\int_{\Omega_k\cap\{u\leq\theta\}} |u(x)+1|^m\,dx+
\int_{B_{kT}\cap \{u\le v_k\}} (v_k(x)+1)^m\,dx\\
\le\;\, & C
\int_{\Omega_k} |v_k(x)+1|^m\,dx+
\int_{B_{kT}} |v_k(x)+1|^m\,dx\\
\le\;\,& C\left(
\int_{\Omega_k\setminus B_{kT}} |v_k(x)+1|^m\,dx+
\int_{B_{kT}} |v_k(x)+1|^m\,dx\right).
\end{split}
\end{equation}
Now we define
\begin{equation}\label{aldAksru4}
\begin{split}&
\alpha_k:=
\int_{B_{kT}} |v_k(x)+1|^m\,dx\\{\mbox{ and }}\;&
\beta_k:=
\int_{\Omega_k\setminus B_{kT}} |v_k(x)+1|^m\,dx=
\int_{(B_{(k+1)T}\setminus B_{kT})\cap\{u>v_k\}} |v_k(x)+1|^m\,dx
.\end{split}\end{equation}
With this, we can rewrite~\eqref{LAKKA:01284} and~\eqref{9qw8e7ryt32938yuwhjsdhfg}
as
\begin{equation}\label{plm9uh8ueygsviudh234}
\begin{split}
&{\mathcal{L}}^n\big(B_{kT}\cap\{u>\theta\}\big)^{\frac{n-1}n}\le
C (\alpha_k+\beta_k)\\
{\mbox{and }}\quad&
\int_{B_{kT}\cap \{u\le\theta\}} |u(x)+1|^m\,dx\le
C (\alpha_k+\beta_k).\end{split}
\end{equation}
We remark that, in view of~\eqref{vk:owue}
and~\eqref{aldAksru4},
\begin{equation}\label{po99df095843ytfjdhsgf}\begin{split}
\alpha_k\;&\le
\int_{B_{kT}} \frac{C}{\big( 1+(k+1)T-|x|
\big)^{\frac{pm}{m-p}}} \,dx\\
&\le C \int_0^{kT} \frac{\rho^{n-1}}{
\big( 1+(k+1)T-\rho
\big)^{\frac{pm}{m-p}}} \,d\rho\\
&\le C k^{n-1} T^{n-1}\int_0^{kT} \frac{d\rho}{
\big( (k+1)T-\rho
\big)^{\frac{pm}{m-p}}} \\
&= C k^{n-1}T^{n-1}\left[ \frac1{T^{\frac{(p-1)m+p}{m-p}}}- \frac{1}{\big( (k+1)T
\big)^{\frac{(p-1)m+p}{m-p}}}\right]\\&\le
C k^{n-1} T^{n-\frac{pm}{m-p}}.
\end{split}
\end{equation}
On the other hand, from~\eqref{aldAksru4}
and the fact that~$v_k(x)\in[-1,1]$ in~$B_{(k+1)T}$,
\begin{eqnarray*}
\beta_k &=&
\int_{(B_{(k+1)T}\setminus B_{kT})\cap\{\theta\ge u>v_k\}} (v_k(x)+1)^m\,dx
+\int_{(B_{(k+1)T}\setminus B_{kT})\cap\{u>\max\{v_k,\theta\}\}} (v_k(x)+1)^m\,dx
\\&\le&
\int_{(B_{(k+1)T}\setminus B_{kT})\cap\{\theta\ge u>v_k\}} (u(x)+1)^m\,dx
+\int_{(B_{(k+1)T}\setminus B_{kT})
\cap\{u>\max\{v_k,\theta\}\}} 2^m\,dx
\\&\le&
\int_{(B_{(k+1)T}\setminus B_{kT})\cap\{\theta\ge u>v_k\}} 
|u(x)+1|^m\,dx
+C{\mathcal{L}}^n \big((B_{(k+1)T}\setminus B_{kT})
\cap\{u>\theta\}\big).
\end{eqnarray*}
Now we insert this information and~\eqref{po99df095843ytfjdhsgf}
into~\eqref{plm9uh8ueygsviudh234}
and we obtain that
\begin{equation}\label{09ouy012ie}
\begin{split}
&{\mathcal{L}}^n\big(B_{kT}\cap\{u>\theta\}
\big)^{\frac{n-1}n}+
\int_{B_{kT}\cap \{u\le\theta\}} |u(x)+1|^m\,dx\\
\le\;\,&
C \left( k^{n-1} T^{n-\frac{pm}{m-p}}+
\int_{(B_{(k+1)T}\setminus B_{kT})\cap\{u\le\theta\}} 
|u(x)+1|^m\,dx
+{\mathcal{L}}^n \big((B_{(k+1)T}\setminus B_{kT})
\cap\{u>\theta\}\big)
\right).
\end{split}\end{equation}
It is now convenient to introduce the following quantities
(somehow reminiscent of ``area'' and ``volume'' terms):
\begin{eqnarray*}
&& {\mathcal{A}}_r:=\int_{B_{rT}\cap\{u\le\theta\}}|u(x)+1|^m\,dx 
\\ {\mbox{and }}
&& {\mathcal{V}}_r:={\mathcal{L}}^n \big(B_{rT}\cap\{u>\theta\}\big)
.\end{eqnarray*}
With this notation, we can rewrite~\eqref{09ouy012ie} as
\begin{equation}\label{023oirgfdoe3u95}
{\mathcal{V}}_k^{\frac{n-1}n}+
{\mathcal{A}}_k\,
\le\,
C \Big( k^{n-1} T^{n-\frac{pm}{m-p}}+
\big( {\mathcal{A}}_{k+1}-{\mathcal{A}}_k\big)
+\big( {\mathcal{V}}_{k+1}-{\mathcal{V}}_k\big)
\Big).
\end{equation}
Now we let~$\varepsilon$ be the structural constant\footnote{We take this opportunity to amend
a typo in~\cite{FV}: the left-hand side of~(61) there should be ``$
{\mathcal{V}}_{k}^{{(n-1)/n}}+{\mathcal{A}}_{k}$'' instead of
 ``${\mathcal{V}}_{k+1}^{{(n-1)/n}}+{\mathcal{A}}_{k+1}$''.}
in Lemma~12 of~\cite{FV}.
Exploiting~\eqref{CONDIZIONE}, for~$T$ large enough we have that~$T^{n-\frac{pm}{m-p}}\le\varepsilon$.
This and~\eqref{023oirgfdoe3u95}, together with~\eqref{INIZIO},
imply the hypotheses of Lemma~12 in~\cite{FV},
which in turn implies that
$$ {\mathcal{A}}_{k}+{\mathcal{V}}_{k}\ge ck^n,$$
for some constant~$c>0$, and this proves~\eqref{FINE}.~\hfill$\Box$

\section{Proof of Corollary~\ref{CORO1}}\label{SCORO1}

The proof of \eqref{n-1} is standard (for instance, one can repeat verbatim the proof
of Lemma~10 in~\cite{FV}).

The proof of~\eqref{LAS} is also a simple consequence of~\eqref{n-1}
and Theorem~\ref{DENS}: we provide full details for the facility of the reader.
We define~$\theta_\star:=\min\{\theta,\vartheta\}$.
Since~$\theta_\star\in(-1,\theta]$, the structural assumptions in~\eqref{SYT:ST1}
and~\eqref{SYT:ST2} hold true with~$\theta$ replaced by~$\theta_\star$.
In addition, since~$\theta_\star\le\vartheta$ we have that~$\{ u>\vartheta\}\subseteq
\{ u>\theta_\star\}$ and thus, by~\eqref{0q9fg09weuty3049yt4},
\begin{equation*}
{\mathcal{L}}^n\big(B_{r_o}\cap\{ u>\theta_\star\}\big)\ge
{\mathcal{L}}^n\big(B_{r_o}\cap\{ u>\vartheta\}\big)\ge c_o
.\end{equation*}
These considerations imply that we can exploit Theorem~\ref{DENS} with~$\theta_\star$
replacing~$\theta$. In this way, from~\eqref{FINE} we deduce that, for large~$r$,
\begin{equation}\label{FINE2}
\int_{B_{r}\cap\{ u\le\theta_\star\}} |u(x)+1|^m\,dx+
{\mathcal{L}}^n\big(B_{r}\cap\{ u>\theta_\star\}\big)\ge \hat c\,r^n
.\end{equation}
On the other hand, by~\eqref{n-1} and~\eqref{N}, for large~$r$ we have that
\begin{eqnarray*}
\frac{C r^{n-1}}{\lambda_\Theta} &\ge& \frac1{\lambda_\Theta}\,{\mathcal{E}}_{B_r}(u)\\
&\ge& \int_{B_r} |u(x)+1|^{{ m }}\chi_{(-\infty,\Theta]}(u(x))\,dx
\\
&\ge& \int_{B_r\cap\{ u\le\theta_\star\} } |u(x)+1|^{{ m }}\,dx+
\int_{B_r\cap\{ \theta_\star < u\le\Theta\}} |u(x)+1|^{{ m }}\,dx\\
\\
&\ge& \int_{B_r\cap\{ u\le\theta_\star\} } |u(x)+1|^{{ m }}\,dx+(1+\theta_\star)^m
{\mathcal{L}}^n\big(B_{r}\cap\{ \theta_\star < u\le\Theta \}\big).
\end{eqnarray*}
That is
$$ \int_{B_r\cap\{ u\le\theta_\star\} } |u(x)+1|^{{ m }}\,dx+
{\mathcal{L}}^n\big(B_{r}\cap\{ \theta_\star < u\le\Theta \}\big)\le \tilde C r^{n-1},$$
with~$\tilde C>0$ also depending on~$\Theta$, $\theta$ and~$\vartheta$.

Hence, recalling~\eqref{FINE2}, for large~$r$ we have that
\begin{eqnarray*}
\hat c\,r^n &\le&
\int_{B_{r}\cap\{ u\le\theta_\star\}} |u(x)+1|^m\,dx
+{\mathcal{L}}^n\big(B_{r}\cap\{ \theta_\star < u\le\Theta \}\big)+
{\mathcal{L}}^n\big(B_{r}\cap\{ u>\Theta\}\big)\\
&\le& \tilde C r^{n-1}+{\mathcal{L}}^n\big(B_{r}\cap\{ u>\Theta\}\big),
\end{eqnarray*}
which implies~\eqref{LAS} when~$r$ is sufficiently large.~\hfill$\Box$

\section{Proof of Corollary~\ref{CORO2}}\label{SCORO2}

The arguments presented for this proof are standard (see e.g.~\cites{MR0445362, CC, MR2126143, MR2139200}).
Since the functional considered here is very general,
we provide the technical details of the proof for the sake of completeness.
We start with the convergence of~$u_\e$
in~$L^1_{\rm loc}(\Omega)$ (this is indeed a general argument,
which does not use the power-like growth from the potential
well, but only a bound from below with a positive 
and continuous function in~$(-1,1)$ vanishing in~$\{1,-1\}$).
We set~$E_o(\tau):=|1-\tau^2|^m$ and
$$ F_o(\tau):=\int_0^\tau \Big( E_o(\sigma)\Big)^{\frac{p-1}p}\,d\sigma.$$
Notice that~$F_o$ is strictly monotone and so we can denote
by~$F_o^{-1}$ the inverse
function of~$F_o$. Let also~$U_\e(x):=F_o(u_\e(x))$.
Notice that
$$ \nabla U_\e(x)= 
\Big( E_o\big(u_\e(x)\big)\Big)^{\frac{p-1}p}\,
\nabla u_\e(x)=
\left(\frac1\e \,E_o\big(u_\e(x)\big)\right)^{\frac{p-1}p}\,
\left(\e^{\frac{p-1}p}\nabla u_\e(x)\right)
$$
and therefore, by Young's Inequality,
$$ |\nabla U_\e(x)| \le
\frac{p-1}{p\e}\,E_o\big(u_\e(x)\big)+
\frac{\e^{p-1}}p\,
|\nabla u_\e(x)|^p.$$
We rewrite this inequality as
$$ c_\star |\nabla U_\e(x)| \le
\frac{\lambda_\star}\e \Big(
\e^p|\nabla u_\e(x)|^p + E_o\big(u_\e(x)\big)
\Big),$$
for a suitable~$c_\star>0$.
{F}rom this, \eqref{SYT:ST1.1} and~\eqref{SUPF}, we conclude that
\begin{equation*}\begin{split}
+\infty\;&>\sup_{\e}  {\mathcal{E}}_\Omega^\e (u_\e)
\\ &\ge \sup_{\e} \frac{\lambda_\star}\e \int_\Omega\Big(
\e^p|\nabla u_\e(x)|^p + E_o\big(u_\e(x)\big)
\Big)\,dx\\
&\ge c_\star\,\sup_{\e}\int_\Omega|\nabla U_\e(x)|\,dx.
\end{split}
\end{equation*}
Accordingly, by the Rellich-Kondrachov Theorem (used here
in~$W^{1,1}(\Omega)$), up to a subsequence
we may suppose that~$U_\e$ converges to some~$U_0$
in~$L^1_{\rm loc}(\Omega)$ and a.e. in~$\Omega$.
Then, we set~$u_0(x):=F^{-1}_o(U_0(x))$ and we have that, a.e.~$x\in\Omega$,
\begin{equation*}
\lim_{\e\searrow0} u_\e(x)=\lim_{\e\searrow0} F^{-1}_o(U_\e(x))=
F^{-1}_o(U_0(x))=u_0(x).
\end{equation*}
Since~$|u_\e|\le1$, this and the Dominated Convergence Theorem
imply that~$u_\e$ converges to~$u_0$ in~$L^1_{\rm loc}(\Omega)$,
as desired.

Now we show that
\begin{equation}\label{ONL}
{\mbox{$u_0$ takes values only in~$\{1,-1\}$.}}\end{equation}
To this end, we use again~\eqref{SYT:ST1.1}
and~\eqref{SUPF}, together with Fatou's Lemma, and we see that
\begin{eqnarray*}
0 &=&\lim_{\e\searrow0} \e {\mathcal{E}}_\Omega^\e(u_\e)\\
&=&\lim_{\e\searrow0}
\int_\Omega E\big(x, \,u_\e(x),\, \e\nabla u_\e(x)\big)\,dx
\\ &\ge& \lambda_\star \lim_{\e\searrow0}
\int_\Omega |1-u_\e^2(x)|^m\,dx
\\ &\ge& \lambda_\star 
\int_\Omega |1-u_0^2(x)|^m\,dx.
\end{eqnarray*}
Consequently, we deduce that~$u_0^2=1$ a.e. in~$\Omega$, which establishes~\eqref{ONL}.

{F}rom~\eqref{ONL}, we have that~$u_\e$ converges to~$\chi_E-\chi_{\R^n\setminus E}$
in~$L^1_{\rm loc}(\Omega)$, with~$E:=\{x\in\Omega{\mbox{ s.t. }}u_0(x)=1\}$.
Now, we prove~\eqref{89:098:98rt}.
To achieve this aim,
we argue by contradiction
and we suppose that there exist~$\Theta\in(0,1)$, $\delta$, $R>0$
and sequences~$\e_k\searrow0$ and~$x_k\in\{ |u_{\e_k}|\le\Theta\}\cap B_R$
with~$B_\delta(x_k)\cap (\partial E)=\varnothing$.
We suppose that~$B_\delta(x_k)\subseteq \R^n\setminus E$ (the case~$B_\delta(x_k)\subset E$
is similar). Then, we have that
\begin{equation}\label{876thjwertyoieruurutututh}
0=\lim_{k\to+\infty}\int_{B_\delta(x_k)} |u_{\e_k}(x)-\chi_E(x)+\chi_{\R^n\setminus E}(x)|\,dx=
\lim_{k\to+\infty}\int_{B_\delta(x_k)} |u_{\e_k}(x)+1|\,dx.
\end{equation}
Now we define~$w_{k}(x):=u_{\e_k}\big(x_k+\e_k x\big)$ and
$$\Omega_k:=\left\{ \frac{x-x_k}{\e_k},\quad x\in\Omega\right\}.$$
{F}rom~\eqref{SUPF11},
\begin{equation}\label{09iuh8yf7tfdx}
\sup_{k\in\N} \| \nabla w_k\|_{L^\infty(\Omega_k)} =
\sup_{k\in\N} \e_k \| \nabla u_{\e_k}\|_{L^\infty(\Omega)}<+\infty.
\end{equation}
We also set
\[
{\mathcal{F}}_{k,\Omega_k} (v):=
\int_{\Omega_k} E\big(\e_k x, \,v(x),\, \nabla v(x)\big)\,dx.
\]
We remark that~$w_k$ is a $Q$-minimum for~${\mathcal{F}}_{k,\Omega_k} $.
Moreover, we have that~$w_k(0)=u_{\e_k}(x_k)\in[-\Theta,\Theta]$.
This and~\eqref{09iuh8yf7tfdx} imply that, fixed~$\Theta'\in(\Theta,1)$,
we have that~$|w_k|\le\Theta'$ in~$B_{r_o}$, for some~$r_o>0$
(in particular, condition~\eqref{0q9fg09weuty3049yt4} is fulfilled here
with~$\vartheta:=-\Theta'$ and~$c_o:={\mathcal{L}}^n(B_{r_o})$).

Therefore, we are in the position of exploiting Corollary~\ref{CORO1}
with the function~$w_k$ and deducing from~\eqref{LAS} (used here 
with~$\Theta:=1/2$)
that,
if~$r$ is sufficiently large,
\begin{equation*}
{\mathcal{L}}^n\big(B_{r}\cap\{ w_k>1/2\}\big)\ge \hat c\,r^n,\end{equation*}
for some~$\hat c>0$.
Hence, scaling back,
\begin{equation*}
{\mathcal{L}}^n\big(B_{\e_k r}(x_k)\cap\{ u_{\e_k}>1/2\}\big)\ge \hat c\,(\e_k r)^n.
\end{equation*}
So, we can choose~$r:=\delta/\e_k$ and plug the latter estimate
into~\eqref{876thjwertyoieruurutututh}. In this way, we obtain that
\begin{equation*}
0\ge\lim_{k\to+\infty}\int_{B_\delta(x_k)\cap\{u_{\e_k}>1/2\}} |u_{\e_k}(x)+1|\,dx\ge
\frac{3}2\,{\mathcal{L}}^n\big(B_{\delta}(x_k)\cap\{ u_{\e_k}>1/2\}\big)
\ge\frac{3\hat c}2\,\delta^n.
\end{equation*}
This is a contradiction, and the proof of~\eqref{89:098:98rt}
is complete.~\hfill$\Box$

\end{document}